\documentclass{aptpub}
\usepackage[applemac]{inputenc}
\usepackage{amssymb}
\usepackage[english]{babel}
 \authornames{}
 \shorttitle{}

\usepackage{amssymb}
\usepackage{amsmath}
\usepackage{fancybox}
\usepackage{graphics}
\usepackage{epsfig}
\usepackage{latexsym}
\usepackage{ae,aecompl}
\usepackage{mathrsfs}
\usepackage{amsfonts,amssymb}

\renewcommand{\E}{{\mathbb{E}}}   
\newcommand{\bt}{{\beta^{\ast}}}

\begin{document}

\title{Partial Match Queries in Two-Dimensional Quadtrees: \\  a Probabilistic Approach.}

\authorone[Département de Mathématiques et Applications, École Normale Supérieure, 45 rue d’Ulm, 75005 Paris, France. e-mail: nicolas.curien@ens.fr ]{Nicolas Curien}
\authortwo[Laboratoire de Probabilités et Modèles Aléatoires,
Université Pierre et Marie Curie, 4 place Jussieu Tour 16-26, 75005 Paris, France. e-mail: adrien.joseph@upmc.fr]{Adrien Joseph}

 \begin{abstract}
 We analyze the mean cost of the partial match queries in random two-dimensional quadtrees. The method is based on fragmentation theory. The convergence is guaranteed by a coupling argument of Markov chains, whereas the value of the limit is computed as the fixed point  of an integral equation. 
 \end{abstract}

 \ \\

\ams{60 F 99}{60 G 18; 60 J 05}

\section{Introduction}
Introduced by  Finkel and Bentley \cite{FB74}, the quadtree structure is a comparison based algorithm designed for retrieving multidimensional data. It is  often studied in computer science because of its numerous applications. The aim of this paper is to study the mean cost of the so-called \emph{partial match queries} in random quadtrees. This problem was first analyzed by Flajolet \emph{et al.}\,\,\cite{FGPR93}.

Let us briefly describe  the discrete model. We choose to focus only on the two-dimensional case. Let $P_{1}, \dots, P_{n}$ be $n$ independent random variables uniformly distributed over $(0,1)^2$. We shall assume that the points have different $x$ and $y$ coordinates, an event that has probability $1$.  We construct iteratively a finite covering of $[0,1]^2$ composed of rectangles with disjoint interiors as follows. The first point $P_{1}$ divides the original square $[0,1]^2$ into four closed quadrants according to the vertical and horizontal positions of $P_{1}$. By induction, a point $P_{k}$ divides the quadrant in which it falls into four quadrants according to its position in this quadrant, see Fig.\,1. Hence the $n$ points $P_{1}, \dots, P_{n}$ give rise to a  covering of $[0,1]^2$ into $3n+1$ closed rectangles with disjoint interiors  that we denote by $\textrm{Quad}(P_{1}, \dots, P_{n})$.

\begin{figure}[!h]
\begin{center}
\includegraphics[width=4cm,height=4cm]{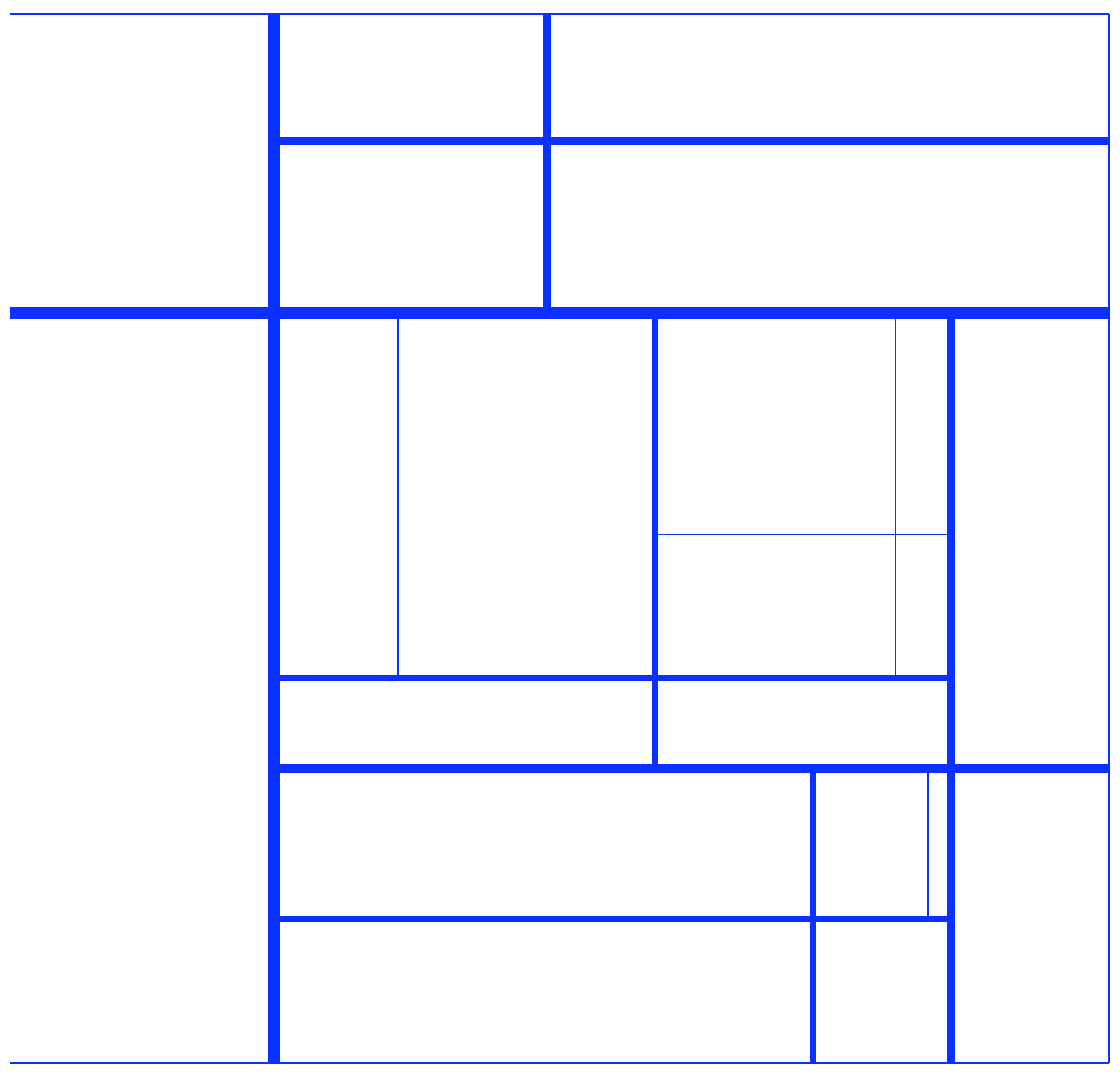} \hspace{2cm}
\includegraphics[width=4cm, height=4cm]{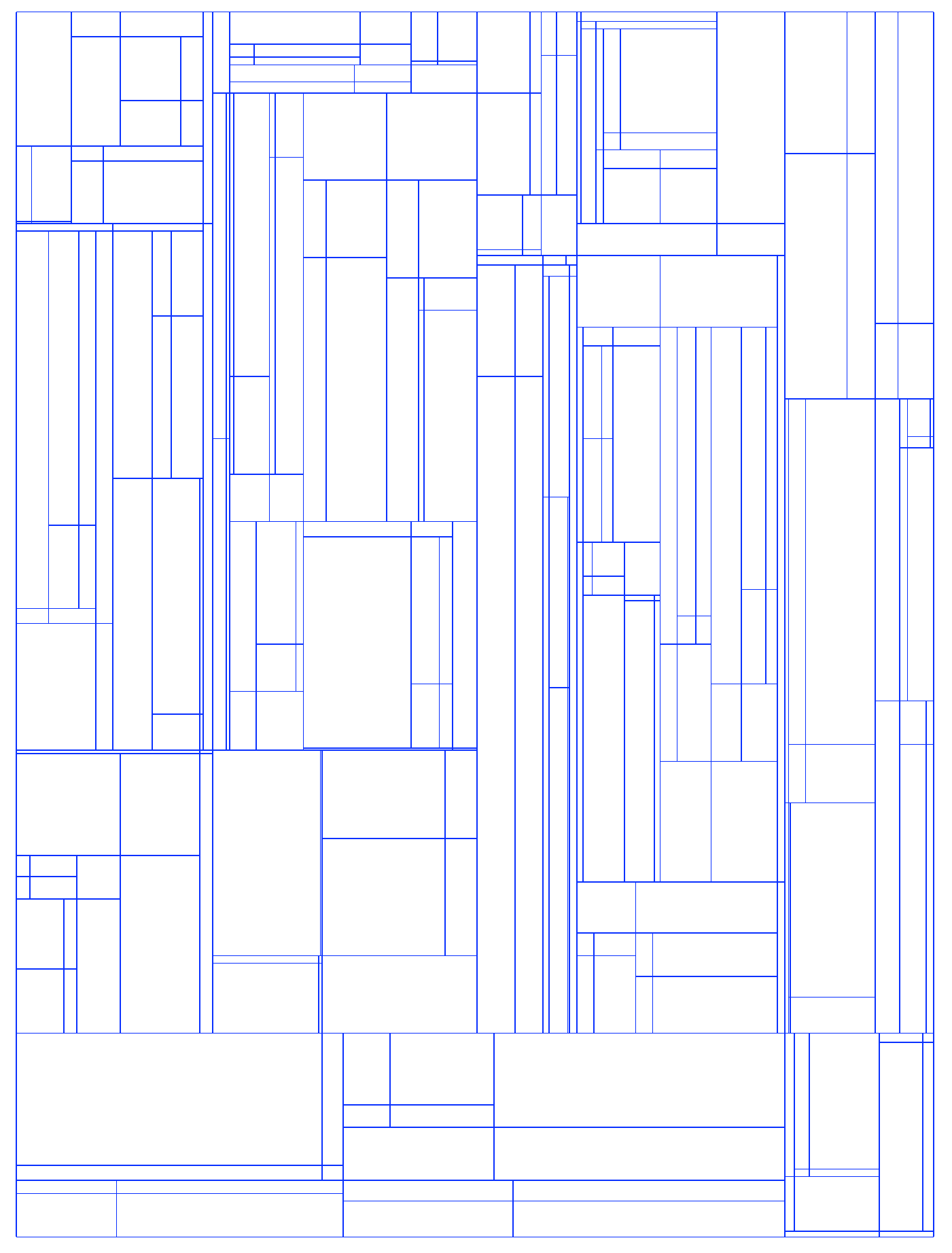}
\caption{Two splittings of $[0,1]^2$ with resp.\,\,8 and 100 points.}
\end{center}
\end{figure}
We are interested in the partial match query. As explained by Flajolet and Sedgewick \cite[Example VII.23.]{FS09},  given $x_0 \in [0,1]$, it  determines the set of points $P_i$, $i \in \{1, \dots, n \}$, with $x$ coordinates equal to $x_0$, regardless of the $y$ coordinates (that set is either empty or a singleton).    Denoting the vertical segment $[(x,0),(x,1)]$ by $S_{x}$, the cost of this partial match query is measured by the number $\mathcal{N}_{n}(x)$ of rectangles of $\textrm{Quad}(P_{1}, \dots, P_{n})$ intersecting $S_{x}$ minus 1 ($\mathcal{N}_{0}(x)=0$ by convention).
Our main result is:

\begin{theorem} \label{mainresul}
For every $x \in [0,1]$, we have the following convergence:
\begin{eqnarray*}
n^{-\beta^{\ast}} \mathbb{E}\big[\mathcal{N}_{n}(x)\big] & \underset{n\to\infty}{\longrightarrow} & K_{0}\big(x(1-x)\big)^{\beta^\ast/2},
\end{eqnarray*}
where  $\displaystyle{\beta^* = \frac{\sqrt{17}-3}{2}}$ and 
$K_{0}=\displaystyle{\frac{\Gamma\left(2\beta^*+2\right) \Gamma(\beta^*+2) }{ 2 \Gamma^3(\beta^*+1) \Gamma^2\left(\frac{\beta^*}{2}+1\right)}}.$
\end{theorem}
Flajolet \emph{et al.}\,\,\cite{FGPR93} obtained the convergence in mean of the cost of partial match queries when $x$ is random with the uniform law on $[0,1]$ and independent of $P_{1}, \dots, P_{n}$. We shall give another proof of this result using \emph{fragmentation theory}. 
 As a by-product of our techniques, we shall prove  in Corollary~\ref{derneircoro} below that when rescaled by $n^{1-\sqrt{2}}$,  $\mathcal{N}_{n}(0)$ converges in $\mathbb{L}^2$ (its convergence in mean was obtained in \cite{FGPR93}).

The paper is organized as follows. Section~\ref{sectiondeux} introduces  the model embedded in continuous-time and presents the first properties. Section~\ref{sectiontrois} is devoted to the link between quadtrees and fragmentation theory. Section~\ref{sectionquatre}, the most technical one, contains the proof of the convergence at a fixed point $x$ without knowing the limit. The identification of the limit is done  in Section~\ref{sectioncinq} using a fixed point argument for integral equation.  \\ 

\textbf{Acknowledgement.} We would like to express our gratitude to Philippe Flajolet who introduced us to the problem of partial match query. We are indebted to Nicolas Broutin and to Ralph Neininger for fruitful discussions. We also deeply thank Jean Bertoin for his careful reading of the first versions of this work. 

\section{Notations and first properties} \label{sectiondeux}
In order to apply probabilistic techniques, we first introduce a continuous-time version of the quadtree: the points $P_{1}, \dots , P_{n}$ are replaced by the arrival points of a Poisson point process over $\mathbb{R}_+ \times [0,1]^2$ with intensity $\textrm{d}t \otimes \textrm{d}x \textrm{d}y$. All the results obtained in this model can easily be translated into results for the discrete-time model.

\subsection{The continuous-time model}

Let $\Pi$ be a Poisson point process on $\mathbb{R}_+ \times [0,1]^2$ with intensity $\textrm{d}t \otimes \textrm{d}x \textrm{d}y$. Let $((\tau_i,x_i,y_i), i\geq 1)$ be the atoms of $\Pi$ ranked in the increasing order of their $\tau$-component. We define a process $(\mathrm{Q}(t))_{t\geq 0}$ with values in finite covering of $[0,1]^2$ by closed rectangles with disjoint interiors as follows. We first introduce the operation $\mathsf{SPLIT}$: for every subset $R$ of $[0,1]^2$ and for every $(x,y) \in [0,1]^2$, 
$$
\mathsf{SPLIT}(R,x,y) =  \big \{ R \cap [0,x] \times [0,y] ,  R \cap [0,x] \times [y,1] , R \cap [x,1] \times [0,y], R \cap [x,1] \times [y,1] \big \}.
$$
In other words, if $R$ is a rectangle with sides parallel to the $x$ and $y$ axes, then $\mathsf{SPLIT}(R,x,y)$ is the set of the four quadrants in $R$ determined by the point $(x,y)$. We may now  recursively define the process $(\mathrm{Q}(t))_{t\geq0}$. Let $\tau_0 = 0$. For every $t \in [0, \tau_1)$, define $\mathrm{Q}(t) = \{ [0,1]^2 \}$, and for every $t \in [\tau_i,\tau_{i+1})$, denoting by $R$ the only element  (if any) of $\mathrm{Q}(\tau_{i-1})$ such that $(x_i, y_i)$ is in the interior of the rectangle $R$,  let
\begin{eqnarray*}
\mathrm{Q}(t) &=& \mathsf{SPLIT}(R,x_i,y_i) \cup \mathrm{Q} \left (\tau_{i-1} \right ) \setminus \{ R\}.
\end{eqnarray*}
Observe that a.s., for every $i \in \mathbb{Z}_+$, there indeed exists a unique  rectangle of $\mathrm{Q}(\tau_i)$ such that $(x_{i+1}, y_{i+1})$ is in its interior, hence the process $(\mathrm{Q}(t))_{t\geq 0}$ is well defined up to an event of zero probability. In the sequel we shall assume that the points of $\Pi$ always fall in the interior of some rectangle of $(\mathrm{Q}(t))_{t\geq 0}$. As explained in the introduction, we are interested in the number of rectangles of $\mathrm{Q}(t)$ intersecting the segment $S_x$, specifically we set:
 \begin{eqnarray*}
N_t(x) & = & \# \big \{ R \in \mathrm{Q}(t) : R \cap S_x \neq \emptyset \big \} -1,
\end{eqnarray*}
so that $N_t(x)=0$ for every $0 \leq t < \tau_1$.  Recalling that $\tau_{n}$ is the arrival time of the $n$-th point of $\Pi$,  $\mathrm{Q}(\tau_{n})$ has the same distribution as the random variable $\mathrm{Quad}(P_{1}, \dots, P_{n})$ of the introduction. In particular,  for every $(n,x) \in \mathbb{N} \times [0,1]$, we have $N_{\tau_{n}}(x) = \mathcal{N}_{n}(x)$ in distribution. 

 \subsection{Main equations}


Let $x \in [0,1]$. We denote by $\mathcal{A}$ the set of words over the alphabet $\{0,1\}$, $$ \mathcal{A} = \bigcup_{n\geq 0} \{ 0,1\}^n,$$ where by convention $\{0,1\}^0= \{ \varnothing\}$. Thus, if $u \in \mathcal{A}$, $u$ is either $\varnothing$ or a finite sequence of $0$ and $1$. If $u$ and $v$ are elements of $\mathcal{A}$ then $uv$ denotes the concatenation of the two words $u$ and $v$. We label the rectangles appearing in $(\mathrm{Q}(t))_{t\geq 0}$ whose intersection with the segment $S_x$ is non-empty by elements of $\mathcal{A}$ according to the following rule. By convention $R^{}_\varnothing(x)$ is the unit square $[0,1]^2$. The first point $(\tau_1,x_1,y_1)$ of $\Pi$ splits $[0,1]^2$ into four rectangles, a.s.\,\,only two of them intersect $S_x$, we denote the bottom rectangle by $R_0(x)$ and the top one by $R_1(x)$. 
Inductively, for every $u \in \mathcal{A}$,  a point of $\Pi$ eventually falls into $R_u(x)$, dividing it into four rectangles. Almost surely, only two of them intersect $S_x$, denote the bottom one by $R_{u0}(x)$ and the top one by $R_{u1}(x)$. For $u \in \mathcal{A}$, we denote  the minimal (resp.\,\,maximal) horizontal coordinate of $R_u(x)$ by $G_u(x)$ (resp.\,\,$D_u(x)$), and define the \emph{place} of $x$ in $R_u(x)$ to be 
$$ X_u(x) = \frac{x-G_u(x)}{D_u(x)-G_u(x)}.$$ If $u \ne \varnothing$, we denote the parent of $u$ by $\overleftarrow{u}$ which is the word $u$ without its last letter. We write $M_u(x)$ for the ratio of the (two-dimensional) Lebesgue measure $\mathrm{Leb}(R_{u}(x))$ of $R_u(x)$ by the measure of $R_{\overleftarrow{u}}(x)$, $$M_u(x) = \frac{\mathrm{Leb}\big(R_u(x)\big)}{\mathrm{Leb}\big(R_{\overleftarrow{u}}(x)\big)}.$$ We also set for all $x\in [0,1]$, $M_{\varnothing}(x)=1$. For $u \in \{0, 1 \}$ and $t\geq 0$, we introduce the ``subquadtree'' $
\mathrm{Q}_{u,x} (t)  =  \left  \{ R \in \mathrm{Q}(t + \tau_1) : R \subset R_u(x) \right \}$. Then, for every $t \geq 0$, one has:
 \begin{eqnarray} \label{premieequt}
N_t(x) & = &\textbf{1}_{t  \geq \tau_1 } + \textbf{1}_{t \geq \tau_1} \sum_{u \in \{0,1\}} \Big ( \# \big\{ R \in \mathrm{Q}_{u,x} (t-\tau_1) : R \cap S_{x} \neq \emptyset \big\} -1 \Big ).
\end{eqnarray}

If $R$ is a rectangle with sides parallel to the $x$ and $y$ axes, we denote by $\Phi_R: \mathbb{R}^2 \to \mathbb{R}^2$ the only affine transformation that maps  the bottom left vertex of $R$ to $(0,0)$,  the bottom right vertex of $R$ to $(1,0)$ and the up left vertex of $R$ to $(0,1)$. It should be plain from properties of Poisson point measures that,  conditionally on $(M_{u}(x),X_{u}(x),R_u(x))$, the process $(\Phi_{R_u(x)}(\mathrm{Q}_{u,x}(t)))_{t\geq 0}$ has the same distribution as the process $(\tilde{\mathrm{Q}}(M_u(x) t))_{t\geq 0}$, where $\tilde{\mathrm{Q}}$ is an independent copy of $\mathrm{Q}$. In particular, conditionally on $(M_{u}(x),X_{u}(x))$,  the number of rectangles in $\mathrm{Q}_{u,x}$ that intersect $S_{x}$ (minus 1), viewed as a process of $t$, has the same distribution as the process $(\tilde{N}_{M_u(x) t}(X_u(x)))_{t\geq 0}$ where $\tilde{N}$ is defined from $\tilde{\mathrm{Q}}$ is the same way as $N$ is defined from $\mathrm{Q}$. Since $M_0(x)$ and $M_1(x)$ have the same distribution, (\ref{premieequt}) yields
 \begin{eqnarray} \label{secoequaf}
\mathbb{E} \left [ N_t(x)  \right ] &=& \mathbb{P}(t \geq \tau_1) + 2 \mathbb{E} \left [ \tilde{N}_{M_{0}(x) (t - \tau_1)}  (X_0(x)  )  \right ],
\end{eqnarray}
with the convention $\tilde{N}_t(x) = 0$ whenever $t < 0$. 
More generally, if we write $\mathfrak{z}_{k} \in \mathcal{A}$ for $\mathfrak{z}_k=0 \dots 0$ repeated $k$ times, then  for every positive integer $k$, 
 \begin{eqnarray} \label{equationfondament}
\mathbb{E} \left [ N_t(x)  \right ] &=& g_k(t) +
2^k \mathbb{E} \left [ \tilde{N}_{M_{\mathfrak{z}_1}(x) \dots M_{\mathfrak{z}_k}(x) t -F_k}  (X_{\mathfrak{z}_k}(x))  \right ],
\end{eqnarray}
where $g_k$ is a function such that $0 \leq g_k \leq 2^k-1$ and $F_k$ is a nonnegative random variable defined by
 \begin{eqnarray*}
F_k &=& \sum_{i=1}^k \tilde{\tau}_i \prod_{j=i}^k M_{\mathfrak{z}_j}(x),
\end{eqnarray*}
with $(\tilde{\tau}_i)_{i\geq1}$  a sequence of independent exponential variables with parameter 1. 

We know compute the joint distribution of $(M_{\textrm{0}}(x), X_0(x))$ which will be of great use throughout this work. If $f$ is a nonnegative measurable function, easy calculations yield
 \begin{eqnarray}
\mathbb{E} \Big[ f \big ( M_{0}(x), X_0(x) \big) \Big ] &=& \int_0^1 \textrm{d} u \int_0^1 \textrm{d} v  \Bigg( \textbf{1}_{x < u} f \left ( uv, \frac{x}{u} \right ) +  \textbf{1}_{x > u} f \left ( (1-u)v, \frac{x-u}{1-u} \right ) \Bigg) \nonumber \\
& = & \int_x^1 \frac{\textrm{d}y}{y} \int_0^{\frac{x}{y}} \textrm{d}m f(m,y) + \int_0^x \frac{\textrm{d}y}{1-y} \int_0^{\frac{1-x}{1-y}} \textrm{d}m f(m,y) \label{loimx}\\
& =& \int_{0}^x \textrm{d}m \int_{x}^1 \frac{\textrm{d}y}{y}f(m,y) + \int_{x}^1 \textrm{d}m\int_{x}^{\frac{x}{m}}\frac{\textrm{d}y}{y}f(m,y) \nonumber \\
 & +& \int_{0}^{1-x} \textrm{d}m\int_{0}^x\frac{\textrm{d}y}{1-y} f(m,y) + \int_{1-x}^{1} \textrm{d}m\int_{1-\frac{1-x}{m}}^x \frac{\textrm{d}y}{1-y} f(m,y) \label{loimx2}.
 \end{eqnarray}
 
 \subsection{Depoissonization}
The following lemma contains a large deviations argument that will enable us to shift results from the continuous-time model to the discrete-time one.

 \begin{lemma} \label{relationnt2}
 For every $\varepsilon>0$, we have
\begin{eqnarray*}
\mathbb{E} \left [ \sup_{x \in [0,1]}  \big | N_{\tau_n}(x) - N_n(x) \big|^2 \mathbf{1}_{\tau_n \notin [n(1-\varepsilon), n(1+\varepsilon)]}  \right ]  &\underset{n \to \infty}{\longrightarrow} & 0.
\end{eqnarray*}
\end{lemma} 
\proof Note that for every $x\in [0,1]$, $t\mapsto N_{t}(x)$ is non-decreasing and that $N_{t}(x)$ is at most the number of points fallen so far:  $N_{t}(x) \leq \max \left \{i \in \mathbb{Z}_+ : \tau_i \leq t \right \}$. In particular $N_{\tau_{n}}(x) \leq n$, thus we have
\begin{eqnarray*}
\sup_{x \in [0,1]}  \big | N_{\tau_n}(x) - N_n(x) \big|^2  \mathbf{1}_{\tau_n > n(1+\varepsilon)}  &\leq& n^2 \mathbf{1}_{\tau_{n} > n(1+\varepsilon)}.
\end{eqnarray*}
A large deviations argument ensures that $n^2 \mathbb{P} ( \tau_n > n(1+\varepsilon) )$ tends to 0 as $n\to\infty$. On the other hand, applying the Cauchy-Schwarz inequality, we obtain
$$
\mathbb{E} \left [ \sup_{x \in [0,1]}  \big | N_{\tau_n}(x) - N_n(x) \big|^2 \mathbf{1}_{\tau_n < n(1+\varepsilon)} \right ]   \leq  \sqrt{ \E \left[ 
(\max \left \{i \in \mathbb{Z}_+ : \tau_i \leq n \right \} )^4 \right]} \sqrt{\mathbb{P} \big( \tau_n < n(1-\varepsilon) \big)}.
$$
As $ 
\E [ (\max \{i \in \mathbb{Z}_+ : \tau_i \leq n  \} )^4 ]
=O(n^4)$, large deviations ensure that the quantity in the right-hand side  tends to 0 as $n \to \infty$. Finally, Lemma~\ref{relationnt2} is proved. 
\endproof

 \section{Particular cases and fragmentation theory} \label{sectiontrois}
 
We give below the definition of a particular case of fragmentation process. For more details, we refer to \cite{Ber06}.
Let $\nu$ be a probability measure on $\{ (s_1, s_2) : s_{1}\geq s_{2} >0 \textrm{ and } s_{1}+s_{2} \leq 1\}$. A self-similar fragmentation $(\mathscr{F}_{t})_{t\geq 0}$ with dislocation measure $\nu$ and index of self-similarity 1 is a Markov process with values in the set $\mathcal{S}^\downarrow = \{ (s_1, s_2, \dots) : s_{1} \geq s_2 \geq \dots  \geq 0 \textrm{ and } \sum_{i} s_{i} \leq 1\}$ describing the evolution of the masses of particles that undergo fragmentation. The process is informally characterized as follows: if at time $t$ we have $\mathscr{F}(t)=(s_{1}(t), s_{2}(t), \dots )$, then for every $i\geq 1$, the $i$-th ``particle'' of mass $s_{i}(t)$ lives an exponential time with parameter $s_{i}(t)$  before splitting into two particles of masses $r_{1}s_{i}(t)$ and $r_{2}s_{i}(t)$, where $(r_{1},r_{2})$ has been sampled from $\nu$ independently of the past and of the other particles. In other words, each particle undergoes a self-similar fragmentation with time rescaled by its mass. In the next section we establish a link between fragmentation theory and the process $N_{t}(U)$, where $U$ is a r.v.\,\,uniformly distributed over $[0,1]$ and independent of $(\mathrm{Q}(t))_{t \geq 0}$. This connection will provide a new proof of a result of \cite{FGPR93} and \cite{CHH03}. See also \cite{CLG10} for another recent  application of fragmentation theory to a combinatorial problem where the exponent $\frac{\sqrt{17}-3}{2}$ appears.

 \subsection{The uniform case}   \label{uniformsection}
   
    We consider here the case where the point $x$ is chosen at random uniformly over $[0,1]$ and independently of $(\mathrm{Q}(t))_{t\geq 0}$.

\begin{proposition} \label{uniuni} 
Let $U$ be a random variable uniformly distributed over $[0,1]$ and independent of the  quadtree $(\mathrm{Q}(t))_{t\geq0}$.  Let $u \in \mathcal{A}$ and denote by $u_{0} = \varnothing,u_{1}, \dots , u_{k}=u$ its  ancestors.  Then $X_{u}(U)$ is uniform over $[0,1]$ and independent of $(M_{u_{1}}(U), \dots , M_{u_{k}}(U))$, which is a sequence of independent random variables all having density $2(1-m)\mathbf{1}_{m \in[0,1]}$.
\end{proposition}
\proof We prove Proposition~\ref{uniuni} by induction on $k$. Let $u \in \mathcal{A}$. Denote by  $u_{0} = \varnothing, u_{1}, \dots, u_{k}=u$ its  ancestors. 
Integrating (\ref{loimx}) for $x\in [0,1]$, we deduce that for every $v \in \{0,1\}$, $X_{v}(U)$ and $M_{v}(U)$ are independent and distributed according to 
 \begin{eqnarray} \label{uny}
\textbf{1}_{u \in [0,1]} \textrm{d}u \otimes\textbf{1}_{m\in [0,1]}2(1-m) \textrm{d}m.
 \end{eqnarray} 
Recalling that $\mathrm{Q}_{u_{1},U}(t)= \{ R \in \mathrm{Q}(t+\tau_{1}): R \subset R_{u_1}(U)\}$, conditionally on $(X_{u_1}(U),M_{u_1}(U))$, the process $\Phi_{R_{u_1}(U)} (\mathrm{Q}_{u_1,U})$  has the same distribution as $(\tilde{\mathrm{Q}}(M_{u_1}(U)t))_{t\geq 0}$, where $\tilde{\mathrm{Q}}$ is an independent copy of $\mathrm{Q}$. Since $X_{u_{1}}(U)$ is uniform over $[0,1]$, we deduce by induction on the subquadtree $\mathrm{Q}_{u_{1},U}$ that $X_{u}(U)$ is uniform over $[0,1]$ and independent of $(M_{u_{2}}(U), \dots , M_{u_{k}}(U))$ which is a sequence of independent r.v.\,\,all having density $2(1-m)\mathbf{1}_{m \in [0,1]}$. Furthermore it is easy to see that 
\begin{eqnarray*}
\mathbb{E}\Big[(X_{u_i}(U),M_{u_i}(U))_{2 \leq i \leq k} \big|(X_{u_{1}}(U), M_{u_{1}}(U)) \Big] &=& \mathbb{E}\Big[(X_{u_i}(U),M_{u_i}(U))_{2 \leq i \leq k} \big|X_{u_{1}}(U)\Big].
\end{eqnarray*} 
Hence by \eqref{uny}, $X_{u}(U)$ is also independent of $M_{u_{1}}(U)$. 
\endproof


Letting $\mathbf{m}(t)= \mathbb{E}[N_t(U)]$, (recall that when $t<0$, $N_{t}(x) =0$ for all $x\in [0,1]$) equation (\ref{secoequaf}) becomes 
\begin{eqnarray}
\label{eq123}\mathbf{m}(t) &=&\mathbb{P}(t \geq \tau_1) + 2 \mathbb{E} \big [ \mathbf{m}(M(t-\tau_1)) \big ],\end{eqnarray} where $M$ is independent of $\tau_1$ and has density $2(1-m)\mathbf{1}_{m\in[0,1]}$.

 \begin{proposition}
 \label{uni} Let $U$ be uniform over $[0,1]$ and independent of $(\mathrm{Q}(t))_{t\geq0}$. We have the following convergence 
\begin{eqnarray*}
\lim_{t \rightarrow \infty} t^{-\beta^*} \E \big[N_t(U)\big] &=& \frac{\Gamma(2(\beta^*+1))}{2\Gamma^3(\beta^*+1)}, \quad \mbox{ where }\beta^*= \displaystyle \frac{\sqrt{17}-3}{2}.
\end{eqnarray*}

\end{proposition}

\proof  We consider an auxiliary fragmentation process $(\mathscr{F}_{t})_{t\geq 0}$ with index of self-similarity 1 and   dislocation probability measure $\nu$ given by 
\begin{eqnarray*} 
\int \nu(ds_1,ds_2) f(s_1,s_2) &=& \E\Big[f\big(M_1(U)\vee M_0(U),M_{1}(U) \wedge M_{0}(U)\big)\Big].
\end{eqnarray*}
In other words, the dislocation measure is given by the law of the decreasing ordering of $\{M_0(U),M_1(U)\}$. More precisely  $(\mathscr{F}_{t})_{t\geq0}$ takes its values in $\mathcal{S}^{\downarrow}$ and satisfies the following equation in distribution which completely characterizes its law:
$$(\mathscr{F}_{t}) \overset{(d)}{=} \bigg (\left (\mathbf{1}_{t < \tau} \right) \dot{+} \left( \mathbf{1}_{t \geq \tau} M_{0}(U)\cdot\mathscr{F}^{(0)}_{M_{0}(U)(t-\tau)}\right)_{t\geq 0} \dot{+}\left( \mathbf{1}_{t \geq \tau}   M_{1}(U)\cdot\mathscr{F}^{(1)}_{M_{1}(U)(t-\tau)}\right)_{t\geq 0} \bigg)^{\downarrow},$$
with $(\mathscr{F}^{(0)}_{t})_{t\geq 0}$ and $(\mathscr{F}^{(1)}_{t})_{t\geq 0}$ two independent copies of $(\mathscr{F}_{t})_{t\geq 0}$  also independent of $(M_{0}(U),M_{1}(U),\tau)$ and $\tau$  an independent exponential variable with parameter $1$. The symbol $\dot{+}$ means concatenation of sequences and $(.)^{\downarrow}$ is the decreasing reordering (and erasing of zeros). 
Then, it is straightforward to see that the expectation of the number $\# \mathscr{F}_{t}$ of fragments of $\mathscr{F}_{t}$ minus $1$ satisfies the same equation as $\mathbb{E}[N_{t}(U)]$, namely letting  $\mathfrak{m}(t)= \mathbb{E}[\#\mathscr{F}_{t}-1]$ for $t \geq 0$, and $\mathfrak{m}(t)=0$ for $t <0$ we have
\begin{eqnarray}
\label{eq124}
\mathfrak{m}(t)&=& \mathbb{P}(t \geq \tau_1) + 2 \mathbb{E} \big[ \mathfrak{m}(M(t-\tau_{1}))\big ], \end{eqnarray}
where $M$ is independent of $\tau_{1}$ and has density $2(1-m)\mathbf{1}_{m\in[0,1]}$. By \eqref{eq123} and \eqref{eq124}, the functions $\mathbf{m}$ and $\mathfrak{m}$ satisfy the same integral equation, 
$$ f(t) = 1-e^{-t} +2 \int_{0}^1 \mathrm{d}m\,2(1-m)\int _{0}^t \mathrm{d}s\,e^{-s}f\big(m(t-s)\big).$$
Differentiating with respect to $t$,
we see that both $\mathbf{m}$ and $\mathfrak{m}$ are solutions of the Cauchy problem for the integro-differential equation 
$$
\left\{ \begin{array}{l}\displaystyle\partial_{t} f(t) = 1-f(t)+ \int_{0}^1 \mathrm{d}m\,2(1-m) f(mt), \\ f(0)=0. \end{array} \right.
$$
Uniqueness of solution of this kind of integro-differential equation is known, see \emph{e.g.}\,\,\cite{IL97}. We deduce that for every $t \geq 0$, $\mathfrak{m}(t) = \mathbf{m}(t)$. We now focus on $\mathfrak{m}(t)$. Following \cite[Section 3]{BG04}, we let for every $\beta >0$, $\psi(\beta) = 1- \int \mathrm\nu(ds_{1},ds_{2})(s_{1}^\beta + s_{2}^\beta)$. An easy calculation yields:
$$
\psi(\beta) = \frac{\beta^2+3\beta-2}{(\beta+1)(\beta+2)}.
$$ 
In particular the Malthusian exponent associated to $\nu$, which is characterized by $\psi(\beta)=0$ (see \cite[Section 1.2.2]{Ber06}),  is $$\beta^*= \frac{\sqrt{17}-3}{2}.$$ Applying \cite[Theorem 1]{BG04}, we get:
$$
\lim_{t \to \infty} t^{-\beta^*}\E[\#\mathscr{F}_t]=
 \frac{\Gamma(1-\beta^\ast)}{\beta^\ast} \frac{4}{2 \beta^\ast + 3}  \prod_{k=1}^\infty \left (1 - \frac{\beta^\ast}{k} \right ) \left (1 - \frac{\beta^\ast}{k+\sqrt{17}} \right ) \left (1 + \frac{\beta^\ast}{k+1} \right ) \left (1 + \frac{\beta^\ast}{k+2} \right ) .
$$
Finally, we use the Weierstrass identity for the gamma function: for every complex number $z \in \mathbb{C} \setminus \mathbb{Z}_{-}$,
$$
\Gamma(z+1) = e^{- \gamma z} \prod_{k=1}^{\infty} \left ( 1+\frac{z}{k} \right )^{-1} e^{z/k},
$$
where $\gamma$ is the Euler–Mascheroni constant. We conclude that
$$
 \lim_{t \to \infty} t^{-\beta^*}\mathbb{E}[N_{t}(U)] =  \frac{4}{\beta^*(2 \beta^\ast + 3)} \frac{\Gamma(\sqrt{17}+1)}{\Gamma(\sqrt{17} - \beta^{\ast} +1)} \frac{1}{\Gamma^2(\beta^{\ast}+2)} \frac{1}{1+\beta^\ast/2}= \frac{\Gamma(2(\beta^*+1))}{2\Gamma^3(\beta^*+1)},
$$
which completes the proof of the proposition.
 \endproof
 
 \begin{remark} \rm
One can derive the following equality in distribution from (\ref{premieequt}):
\begin{eqnarray*}
N_t(U) &\overset{(d)}{=} &\mathbf{1}_{\tau_1 \leq t} + N^{(0)}_{M_0(U)(t-\tau_1)}\big(X_0(U)\big) +N^{(1)}_{M_1(U)(t-\tau_1)}\big(X_1(U)\big),
\end{eqnarray*} 
where $(N^{(0)}_{t})_{t\geq 0}$ and $(N^{(1)}_{t})_{t\geq 0}$ are independent copies of the process $(N_{t})_{t\geq 0}$. We have already noticed that $X_0(U)$ and $X_1(U)$ are also uniform and independent of $(N^{(0)}_{t})_{t\geq 0}$, of $(N^{(1)}_{t})_{t\geq 0}$ and of $(M_{0}(U),M_{1}(U))$. If $X_{0}(U)$ and $X_{1}(U)$ were independent, then $N_{t}(U)$ would satisfy the same distributional equation as $(\#\mathscr{F}_{t}-1)_{t \geq 0}$. However, this is not the case since we have $X_0(U)=X_1(U)$. This explains why we had to work with expectations.
\end{remark}

\begin{corollary}[\cite{FGPR93}, \cite{CHH03}] \label{remarkcompute} \rm We have
\begin{eqnarray*}
\lim_{n \rightarrow \infty} n^{-\beta^*} \E\big[\mathcal{N}_{n}(U)\big] & = & \frac{\Gamma(2(\beta^*+1))}{2\Gamma^3(\beta^*+1)}.
\end{eqnarray*}
\end{corollary}
\proof This is a straightforward application of Lemma~\ref{relationnt2} and  Proposition~\ref{uni}. \endproof


\subsection{Case $x=0$}
As a further example of the connection with fragmentation theory, we derive asymptotics properties for $N_{t}(0)$. In this case, the sequence of the areas of the rectangles crossed by $S_{0}$ \emph{is} a fragmentation process, enabling us to state a  convergence of $N_{t}(0)$, once rescaled, in $\mathbb{L}^2$. A convergence in mean has already been obtained in \cite[Theorem 6]{FGPR93} and \cite{FLLS95}.
\begin{theorem} \label{lastthro}
The random variable 
$$ \mathfrak{M}_t  =  \sum_{u \in \mathcal{A}} \mathrm{Leb} \big( R_u(0) \big)^{\sqrt{2}-1}  \mathbf{1}_{R_u(0) \in \mathrm{Q}(t)}, \quad t \geq 0,$$
is a uniformly integrable martingale which converges almost surely to $\mathfrak{M}_{\infty}$ as $t\to\infty$. The distribution of $\mathfrak{M}_{\infty}$ is   characterized by
\begin{eqnarray}
\label{eqloiloi}
\mathbb{E}[\mathfrak{M}_\infty]=1 \quad \textrm{and} \quad \mathfrak{M}_{\infty} \overset{(d)}{=} M_0(0)^{\sqrt{2}-1}\mathfrak{M}^{(0)}_{\infty} +   M_1(0)^{\sqrt{2}-1}\mathfrak{M}^{(1)}_{\infty},
\end{eqnarray}
where $\mathfrak{M}^{(0)}_{\infty}$ and $\mathfrak{M}^{(1)}_{\infty}$ are two independent copies of  $\mathfrak{M}_{\infty}$ also independent of  $(M_0(0),M_1(0))$. Furthermore, we have the following convergence in $\mathbb{L}^2$:
\begin{eqnarray*}
t^{1-\sqrt{2}} N_{t}(0) &\underset{t \to \infty} { \longrightarrow}& \frac{\Gamma(2 \sqrt{2})}{\sqrt{2}\Gamma^3(\sqrt{2})} \mathfrak{M}_\infty.
\end{eqnarray*}
 \end{theorem}
\proof It is easy to check from properties of Poisson measures that the rearrangement in decreasing order of the masses of the rectangles living at time $t$ and intersecting $S_{0}$,
$$ \Big( \mathrm{Leb}\big(R_{u}(0)\big)\mathbf{1}_{R_{u}(0) \in \mathrm{Q}(t)} \Big)_{t \geq 0}^{\downarrow},$$ is a self-similar fragmentation with index 1 and  dislocation probability measure given by the decreasing ordering of $\{M_0(0),M_1(0)\}$. As in the proof of Proposition~\ref{uni}, we introduce for every $\beta>0$, $\Psi(\beta) = 1 - \mathbb{E}[M_{0}(0)^\beta + M_{1}(0)^\beta]$, which is easily computed:
$$ \Psi(\beta) = \frac{(\beta+1)^2-2}{(\beta+1)^2}.$$
Thus the Malthusian exponent $p^\ast$ of this fragmentation satisfying $\Psi(p^*)=0$ is $$p^*=\sqrt{2}-1.$$  The first two points of the theorem follow from classical results of fragmentation theory, see \cite[Theorem 1.1]{Ber06}. We refer to \cite{Liu97} for the characterization of the law of $\mathfrak{M}_{\infty}$ \textit{via} the distributional equation \eqref{eqloiloi} and to \cite{Liu01} for some of its properties. The last point comes from \cite[Corollary 6]{BG04} and the Weierstrass identity for the gamma function used in a similar manner as in the proof of Proposition~\ref{uni}.  
\endproof


\begin{corollary}  \label{derneircoro}
We have the following convergence in $\mathbb{L}^2$:
\begin{eqnarray*}
n^{1-\sqrt{2}} N_{\tau_n}(0) &\underset{n \to \infty} { \longrightarrow}& \frac{\Gamma(2 \sqrt{2})}{\sqrt{2}\Gamma^3(\sqrt{2})} \mathfrak{M}_\infty.
\end{eqnarray*}
\end{corollary}

\proof   This proposition easily derives from Lemma~\ref{relationnt2} and Theorem~\ref{lastthro}. 
 \endproof

\begin{remark} \rm
Observe that Corollary~\ref{derneircoro} implies the following convergence in distribution:
\begin{eqnarray*}
n^{1-\sqrt{2}} \mathcal{N}_{n}(0) &\underset{n \to \infty} { \longrightarrow}& \frac{\Gamma(2 \sqrt{2})}{\sqrt{2}\Gamma^3(\sqrt{2})} \mathfrak{M}_\infty.
\end{eqnarray*}
 \end{remark}

\begin{remark} \rm
It is worthwhile to notice that the behavior of the cost of the partial match query in the case $x=0$ is drastically different from its behavior in the case when $x$ is uniform or $x$ is fixed in $(0,1)$ (see Theorem~\ref{mainresul} and Proposition~\ref{uni}).
\end{remark}

\subsection{An \textit{a priori} uniform bound}

This section is devoted to the proof of an \textit{a priori} uniform bound on $s^{-\beta^*}\mathbb{E}[N_t(x)]$ over $(x,s) \in (0,1) \times (0, \infty)$ that will be useful in many places.
 \begin{lemma} \label{lemmebornit}
 There exists $C < \infty$ such that
 \begin{eqnarray} \label{bornitude}
 \sup_{x \in (0,1)} \sup_{s > 0}  \ \mathbb{E} \Big [s^{-\beta^{\ast}} N_s(x)\Big] & \leq & C.
\end{eqnarray}
\end{lemma}
\proof As a warmup, we start by proving that there exists $C_1 < \infty$ such that for every $x \in (0,1)$,
 \begin{eqnarray} \label{lastequato} 
\sup_{s > 0}  \ \mathbb{E} \Big [s^{-\beta^{\ast}} N_s(x)\Big] & \leq & \frac{C_{1}}{x \wedge (1-x)}.
\end{eqnarray}
Combining (\ref{secoequaf}) with the densities computed in (\ref{loimx}), we deduce that for every $x \in (0,1)$
\begin{eqnarray}
t^{-\beta^\ast}  \mathbb{E}[N_t(x)] & = &t^{-\beta^\ast} \mathbb{P}(t \geq \tau_1)  +2 \left(\int_x^1  \frac{\textrm{d}y}{y} \int_0^{\frac{x}{y}} \textrm{d}m \mathbb{E} \left [ t^{-\bt}N_{m(t-\tau_1)}(y) \right ] \right. \nonumber \\ 
& & + \left.\int_0^x \frac{\textrm{d}y}{1-y} \int_0^{\frac{1-x}{1-y}} \textrm{d}m \mathbb{E} \left [ t^{-\bt}N_{m(t-\tau_1)}(y) \right ]  \right). \label{pointeight}
\end{eqnarray}
By monotony of $t\mapsto N_{t}(x)$ we have $\mathbb{E} \left [ t^{-\bt}N_{m(t-\tau_1)}(y) \right ] \leq \mathbb{E} \left [ t^{-\bt}N_{t}(y) \right ]$. Furthermore, recalling that $\beta^* <1$, there exists a constant $C'$ such that for every $t > 0$, $t^{-\beta^*}\mathbb{P}(t \geq \tau_{1}) \leq C'$. Hence
\begin{eqnarray*}
t^{-\beta^\ast}  \mathbb{E}[N_t(x)] & \leq & C' + 2 \left( \int_{x}^1 \frac{x \textrm{d}y}{y^2} \mathbb{E} \left [ t^{-\bt}N_{t}(y) \right ] + \int_{0}^x  \frac{(1-x)\textrm{d}y}{(1-y)^2} \mathbb{E} \left [ t^{-\bt}N_{t}(y) \right ] \right) \\
& \leq & C' + \frac{2}{x \wedge (1-x)} \int_{0}^1 \textrm{d}y \mathbb{E} \left [ t^{-\bt}N_{t}(y) \right ] \\
&=&  C' + \frac{2}{x \wedge (1-x)} \mathbb{E}\big[ t^{-\beta^*}N_{t}(U)\big].\end{eqnarray*}
It has been shown in Proposition~\ref{uni} that  $\mathbb{E} \left [ t^{-\bt}N_{t}(U) \right ]$ has a finite limit as $t \to \infty$, and for every $t > 0$, $\mathbb{E} \left [ N_{t}(U) \right ] \leq t$. Thus the quantity $\mathbb{E} \left [ t^{-\bt}N_{t}(U) \right ]$ is bounded  over $(0, \infty)$. The inequality \eqref{lastequato} follows from these considerations.

Introducing $S(x) = \sup_{s>0}s^{-\beta^*}\mathbb{E}[N_{s}(x)]$ for every $x\in [0,1]$, we have just shown that $S(x) \leq C_{1}(x\wedge(1-x))^{-1}$. Using \eqref{pointeight}, we have for every $x \in (1/2,1)$:
\begin{eqnarray}
S(x) & = &\sup_{t >0} \left\{t^{-\beta^\ast} \mathbb{P}(t \geq \tau_1)  +2 \left(\int_x^1  \frac{\textrm{d}y}{y} \int_0^{\frac{x}{y}} \textrm{d}m \mathbb{E} \left [ t^{-\bt}N_{m(t-\tau_1)}(y) \right ] \right.\right. \nonumber \\ & & + \left.\left.\int_0^x \frac{\textrm{d}y}{1-y} \int_0^{\frac{1-x}{1-y}} \textrm{d}m \mathbb{E} \left [ t^{-\bt}N_{m(t-\tau_1)}(y) \right ]  \right)\right\} \nonumber \\
 & \leq & C' + 2 \sup_{t>0} \left\{ \int_x^1  \frac{\textrm{d}y}{y} \int_0^{1} \textrm{d}m \mathbb{E} \left [ t^{-\bt}N_{t}(y) \right ] + \int_0^{1/2} \frac{\textrm{d}y}{1-y} \int_0^{1} \textrm{d}m \mathbb{E} \left [ t^{-\bt}N_{t}(y) \right ] \right \} \nonumber \\ 
 & & +  2 \sup_{t>0} 
  \int_{1/2}^x \frac{\textrm{d}y}{1-y} \int_0^{\frac{1-x}{1-y}} \textrm{d}m m^{\beta^\ast} \mathbb{E} \left [ (mt)^{-\bt}N_{mt}(y) \right ] \nonumber \\
  & \leq & C' + 8 \sup_{t>0} \int_0^1  \textrm{d}y   \mathbb{E} \left [ t^{-\bt}N_{t}(y) \right ]  +  2 
  \int_{1/2}^x \frac{\textrm{d}y}{1-y} \int_0^{\frac{1-x}{1-y}} \textrm{d}m\,m^{\beta^\ast} S(y) \nonumber \\
 & \leq & C_{2} +\frac{2}{\beta^*+1} (1-x)^{\beta^\ast +1}   \int_{1/2}^x \textrm{d}y \frac{1}{(1-y)^{\beta^\ast + 2}} S(y). \label{youpi}
 \end{eqnarray}
  Let us show that this implies that for every $x \in (0,1)$, $S(x) \leq 100 C_2$. Arguing by contradiction, suppose that there exists $a \in (1/2,1)$ such that $S(a) > 100 C_2$. Let $S = \sup_{x \in [1/2, a]} S(x)$. By \eqref{lastequato}, $S$ is finite; there exists $b \in [1/2, a]$ such that $S(b) \geq 0.9 S$. In particular, $S(b) \geq 0.9 \sup_{x \in [1/2, b]} S(x)$ and $S(b) > 90 C_2$. Applying \eqref{youpi} at $b$, we get 
 \begin{eqnarray*} 
 S(b) &\leq& 90^{-1} S(b) + \frac{2}{\beta^*+1} (1-b)^{\beta^\ast +1}   \int_{1/2}^b \textrm{d}y \frac{1}{(1-y)^{\beta^\ast + 2}} 0.9^{-1} S(b) \\
 &\leq& 90^{-1} S(b) + \frac{2 \cdot 0.9^{-1}} {(\beta^*+1)^2} S(b),
 \end{eqnarray*}
leading to a contradiction since $(\beta^* + 1)^2 >\frac{2 \cdot 0.9^{-1}}{1-90^{-1}}$. Finally, $S(x) \leq 100 C_2$ for every $x \in (0,1)$.
 \endproof

\section{The convergence at fixed $x \in (0,1)$} \label{sectionquatre}

We prove in this section that when $x \in [0,1]$ is fixed, $t^{-\beta^\ast} \mathbb{E} [N_t(x)]$ admits a finite limit as $t \to \infty$. The results of the preceding section do not directly apply since the place $X_{0}(x)$ of $x$ in the rectangle $R_{0}(x)$ highly depends on $x$. Recall notation $\mathfrak{z}_{k}$ for the word composed of $k$ zeros $0 \dots 0 \in \mathcal{A}$. The guiding idea is that the splittings tend to make $X_{\mathfrak{z}_{k}}(x)$ uniform and independent of $M_{\mathfrak{z}_k}(x)$.

\subsection{A key Markov chain}   \label{markovchain}


Fix $x \in (0,1)$. To simplify notation, for every $k\geq 1$, we write $X_{k}$ for $X_{\mathfrak{z}_{k}}(x)$ and $M_{k}$ for $M_{\mathfrak{z}_{k}}(x)$. We shall focus on the process $(X_{k},M_{k})_{k\geq 0}$, which is obviously a homogeneous Markov chain starting from $(x,1)$ whose transition probability is given by \eqref{loimx} or \eqref{loimx2}. Let $k\geq1$. We denote by $\mathcal{F}_{k}$ the filtration generated by $(X_{i},M_{i})_{1 \leq i \leq k}$. It is easy to see that the transition probability only depends on $X_k$, that is 
\begin{eqnarray*}
\mathbb{E}\big[(X_{k+i}, M_{k+i})_{i\geq 1} | \mathcal{F}_{k}\big] = \mathbb{E}\big[ (X_{k+i}, M_{k+i})_{i\geq 1} | X_{k}\big].
\end{eqnarray*}

\begin{prop} \label{mc} Fix $x \in (0,1)$. There exists a coupling of the chain $(X_{k}, M_{k})_{k\geq 0}$ with a random time $T \in \mathbb{Z}_{+}$  such that for any $k \geq 0$, conditionally on $\{ T \leq k\}$, the r.v.\,\,$X_{k}$ is uniformly distributed over $[0,1]$,  independent of $(M_{i})_{1 \leq i \leq k}$ and of $T$. 
Furthermore, we have 
 \begin{eqnarray*}
 \mathbb{E}\left[1.15^T \right] < +\infty. 
 \end{eqnarray*}
\end{prop}
\proof 
 For any $k\geq 1$ we consider the event
  $$
  E_k= \left \{ M_{k}< X_{k-1} \wedge (1-X_{k-1}) \right \}.
  $$
  Using the explicit densities  (\ref{loimx}) and (\ref{loimx2}),  one sees that conditionally on $\mathcal{F}_{k-1}$ and on the event $E_k$  of probability $-(X_{k-1}\wedge(1-X_{k-1}))\ln(X_{k-1}(1-X_{k-1})),$  the conditional distribution of $X_k$ is $$\frac{1}{-\ln\big(X_{k-1}(1-X_{k-1})\big)}\left(\frac{1}{1-y}\textbf{1}_{y \in (0,X_{k-1})}+ \frac{1}{y} \textbf{1}_{y \in (X_{k-1},1)}\right) \textrm{d}y.$$ In particular, conditionally on $E_k$ and $\mathcal{F}_{k-1}$,  the variable $X_{k}$ is independent of $M_{k}$ and  has a density bounded from below by $- 1 / \ln(X_{k-1}(1-X_{k-1}))$. Thus, we can construct simultaneously with $(X_{k},M_{k})_{k\geq 0}$ a sequence of random variables $(B_{k})_{k \geq 0} \in \{0,1\}^{\mathbb{Z}_+}$ as follows. Suppose that we have constructed $(X_{i}, M_{i}, B_{i})_{0\leq i\leq k-1}$. Then independently of $\mathcal{F}_{k-1}$, 
  toss a Bernoulli variable of parameter $-(X_{k-1}\wedge(1-X_{k-1}))\ln(X_{k-1}(1-X_{k-1})).$ If $0$ comes out, we consider that we are on the event $E_k^c,$ then put $B_{k}=0$ and sample $(X_{k}, M_{k})$ with the conditional distribution on $E_k^c$ and $\mathcal{F}_{k-1}$. If $1$ comes out,  we consider that we are on the event $E_k$ and we proceed to the following. 
   \begin{enumerate}
 \item First sample $M_{k}$ from its distribution conditionally on $E_k$ and $\mathcal{F}_{k-1}$.
 \item Then independently of $M_{k}$, toss a Bernoulli variable $B_{k}$ of parameter $-1/\ln(X_{k-1}(1-X_{k-1}))$. If $B_{k}=1$, sample $X_{k}$ uniformly from $[0,1]$ and independently of $(M_{1}, \dots ,M_{k})$. Otherwise, sample $X_{k}$ with density 
 $$\frac{1}{-\ln\big(X_{k-1}(1-X_{k-1})\big)-1}\left(  \left (  \frac{1}{1-y} -1 \right ) \textbf{1}_{y \in (0,X_{k-1})}+ \left (  \frac{1}{y} - 1 \right ) \textbf{1}_{y \in (X_{k-1},1)}\right) \textrm{d}y,$$
 independently of $(M_{1}, \dots ,M_{k})$.
 \end{enumerate}
 The device provides us with a Markov chain $(X_{k}, M_{k},B_{k})_{k\geq 0}$ such that the first two coordinates have the law of the process introduced before Proposition~\ref{mc}. We then let $$T= \inf\{ k \geq 0, B_{k}=1 \}.$$ By definition of $T$, the random variable $X_{T}$ is sampled uniformly over $[0,1]$ and independently of $(M_{1}, \dots , M_{T})$.  We deduce that  the process $(X_{T+i}, M_{T+i})_{i\geq 1}$ has the same distribution as the process $(X_{\mathfrak{z}_{k}}(U), M_{\mathfrak{z}_{k}}(U))_{k \geq 1}$ defined in Proposition~\ref{uniuni}, hence an easy adaptation of Proposition~\ref{uniuni} shows that for every positive integer $i$, $X_{T+i}$ is uniformly distributed over $[0,1]$ independent of $(M_{1}, \dots , M_{T+i})$ and of $T$. This proves the first part of Proposition~\ref{mc}. 
 
 For the second part, we need to evaluate the tail of the random time $T$. 
We introduce the following variation. Let $(\hat{X}_{k})_{k\geq0}$ be a Markov chain with space state $(0,1)\cup\{\partial\}$, where $\partial$ is a cemetery point. Informally, this chain is the chain $(X_{k})$ until we reach the time $T$, then it is killed and sent to the cemetery point. 
Thanks to the calculation presented at the beginning of the proof, it should be clear that given $X_{k-1}$ and conditionally on $\{ T \geq k-1\}$, the probability of the event $\{T=k\}$ is 
$X_{k-1} \wedge (1-X_{k-1}).$
Thus the transition probability for the chain $(\hat{X}_{k})$ is defined as follows: for every $x \in (0,1)$,
 \begin{eqnarray*}
p(x, \textrm{d} y)&=& x \wedge (1-x) \delta_{\partial} + \left ( \frac{1-x}{(1-y)^2} \textbf{1}_{y \in (0,x)} + \frac{x}{y^2} \textbf{1}_{y \in (x,1)}   - x \wedge (1-x) \right ) \textrm{d} y,
 \end{eqnarray*}
and
 $p(\partial, \textrm{d} y)=\delta_{\partial}.$
By construction of this chain, the stopping time $\hat{T}= \inf\{k \geq 1: \hat{X}_{k}=\partial\}$ has the same distribution as $T$. In order to estimate $\hat{T}$, we define the following potential function $V : (0,1) \cup \{ \partial \} \to [1, \infty]$:
 \begin{eqnarray*}
V(x) & = &  \textbf{1}_{x = \partial} + \frac{10}{\sqrt{x}} \textbf{1}_{x \in (0,1/2)} +  \frac{10}{\sqrt{1-x}} \textbf{1}_{x \in [1/2,1)}.
 \end{eqnarray*}
 Then one can show that for every $x \in (0,1) \cup \{ \partial \} $,
  \begin{eqnarray*}
\int p(x, \textrm{d} y )V(y) & \leq  & 0.85 V(x) + \textbf{1}_{\{ \partial \}}(x),
 \end{eqnarray*}
so that \cite[Theorem 15.2.5]{MT09} may be applied: there exists $\varepsilon > 0$ such that for all $x \in (0,1)$,
  \begin{eqnarray*}
\mathbb{E} \left [ \sum_{k=0}^{\hat{T}-1} V \left (\hat{X}_k \right ) 1.15^k \right ] & \leq & \varepsilon^{-1} 1.15^{-1} V(x),
 \end{eqnarray*}
 from which we deduce that 
 $$
 \mathbb{E} \left  [1.15^{\hat{T}}  \right ]  <  \infty
 $$
  (note that the last quantity is not uniformly bounded for $x\in (0,1)$). This completes the proof of Proposition~\ref{mc}. \endproof

  
 In the remaining part of this section, $x$ is fixed in $(0,1)$. Coming back to (\ref{equationfondament}) and  writing $\overline{M}_{k}=M_{1} 
 M_{2} \dots M_{k}$ for the Lebesgue measure of $R_{\mathfrak{z}_{k}}(x)$, we have
 \begin{eqnarray} \label{keyequation}
t^{- \beta^\ast}\mathbb{E} \left [ N_t(x)  \right ] = t^{- \beta^\ast} \Big ( g_k(t) +
 2^k \mathbb{E} \left [ \tilde{N}_{\overline{M}_k t -F_k}  (X_k) \textbf{1}_{T > k} \right ] + 
2^k \mathbb{E} \left [ \tilde{N}_{\overline{M}_k t -F_k}  (X_k) \textbf{1}_{T \leq k} \right ] \Big).
\end{eqnarray}
We shall treat separately the last two terms of (\ref{keyequation}). 

\subsection{Study of  $t^{- \beta^\ast} 2^k \mathbb{E} [ \tilde{N}_{\overline{M}_k t -F_k}  (X_k) \textbf{1}_{T > k} ] $}

We shall see that $t^{- \beta^\ast} 2^k \mathbb{E} [ \tilde{N}_{\overline{M}_k t -F_k}  (X_k) \textbf{1}_{T > k} ] $ is arbitrarily  small uniformly in $t$ provided that the integer $k$ is chosen large enough. Observe
 \begin{eqnarray*}
t^{- \beta^\ast} 2^k \mathbb{E} \big [ \tilde{N}_{\overline{M}_k t -F_k}  (X_k) \textbf{1}_{T > k} \big ] &\leq &t^{- \beta^\ast} 2^k \mathbb{E} \big [ \tilde{N}_{\overline{M}_{k} t}  (X_k) \textbf{1}_{T > k} \big ] \\
& = &2^k \mathbb{E} \left [ \overline{M}_{k}^{\beta^\ast} (\overline{M}_{k} t)^{- \beta^\ast}  \tilde{N}_{\overline{M}_{k} t}  (X_k) \textbf{1}_{T > k} \right ] \\
& = & 2^k \mathbb{E} \left [  \overline{M}_{k}^{\beta^\ast} \textbf{1}_{T > k}   \mathbb{E} \left [ \left . (\overline{M}_{k} t)^{- \beta^\ast}  \tilde{N}_{\overline{M}_{k} t}  (X_k)  \right | \sigma(\overline{M}_{k}, X_k, T)  \right ] \right ].
\end{eqnarray*}
Letting $\phi$ be the map $(s,u) \mapsto \mathbb{E} [ s^{-\beta^\ast} N_s(u)]$, we have:
\begin{eqnarray*}
t^{- \beta^\ast} 2^k \mathbb{E} \big [ \tilde{N}_{\overline{M}_k t -F_k}  (X_k) \textbf{1}_{T > k} \big ] &\leq & 2^k \mathbb{E} \left [  \overline{M}_{k}^{\beta^\ast} \textbf{1}_{T > k}   \phi(\overline{M}_{k} t, X_k) \right ].
\end{eqnarray*}
Thanks to \eqref{bornitude}, $\phi \leq C$, so that the quantity in the last display is at most
$C 2^k \mathbb{E} \left [  \overline{M}_k^{\beta^\ast} \textbf{1}_{T > k} \right ].$
Hölder's inequality yields for every $p>1$
\begin{eqnarray*}
C 2^k \mathbb{E} \left [  \overline{M}_k^{\beta^\ast} \textbf{1}_{T > k} \right ] &\leq & C 2^k \mathbb{E} \left [  \overline{M}_k^{\beta^\ast p} \right ]^{1/p} \mathbb{E} \left [ \textbf{1}_{T > k} \right ]^{1-1/p}.
\end{eqnarray*}
The last term is easily treated, by Markov's inequality we have 
 $\mathbb{E} \left [ \textbf{1}_{T > k} \right ]  \leq 1.15^{-k} \E[1.15^T]$. Concerning $\E[\overline{M}_{k}^{\beta^*p}]$ we have
\begin{eqnarray*}
\mathbb{E} \left [  \overline{M}_k^{\beta^\ast p} \right ]& \leq & \mathbb{E} \left [  M_{\mathfrak{z}_{2}} (x)^{\beta^\ast p} \dots M_{\mathfrak{z}_{k}}(x)^{\beta^\ast p} \right ] \\
& = & \int_0^1 f^{(x)}(y) \textrm{d} y  \mathbb{E} \left [ M_{\mathfrak{z}_{1}} (y)^{\beta^\ast p} \dots M_{\mathfrak{z}_{k-1}} (y)^{\beta^*p}  \right ],
\end{eqnarray*}
where $f^{(x)}$ is the density of $X_1$ under $\mathbb{P}$. It is easy to see from (\ref{loimx}) that $f^{(x)}$ is bounded from above by $(x \wedge (1-x))^{-1}$. Hence
\begin{eqnarray*}
\mathbb{E} \left [  \overline{M}_k^{\beta^\ast p} \right ]& \leq & \frac{1}{x \wedge (1-x)} \int_0^1  \textrm{d} y  \mathbb{E} \left [ \overline{M}_{k-1}(y)^{\beta^\ast p}  \right ].
\end{eqnarray*}
Recall from Proposition~\ref{uniuni} that when $x=U$ is uniform over $[0,1]$ and independent of $(\mathrm{Q}(t))_{t\geq 0}$,  then $M_{\mathfrak{z}_{1}}(U), \dots , M_{\mathfrak{z}_{k}}(U)$ are independent and distributed according to $\mathbf{1}_{m \in [0,1]}2(1-m) \textrm{d}m$. In particular
$$
\E\Big[M_{0}(U)^{\beta^*p}\Big]  = \frac{2}{(\beta^*p+1)(\beta^*p+2)}
$$
and thus 
$$
\int_0^1  \textrm{d} y  \mathbb{E} \left [ \overline{M}_{k-1}(y)^{\beta^\ast p}  \right ] =\left(\frac{2}{(\beta^*p+1)(\beta^*p+2)}\right)^{k-1}.
$$
 Gathering all these estimates, we obtain
\begin{eqnarray*}
 & & t^{- \beta^\ast} 2^k \mathbb{E} \big [ N_{\overline{M}_k t -F_k}  (X_k) \textbf{1}_{T > k} \big ]   
\\
 & \leq & C 2^k \left(   \frac{1}{x \wedge (1-x)}      \right )^{1/p}    \left ( \frac{2}{(\beta^\ast p + 1)(\beta^\ast p + 2)} \right )^{(k-1)/p}   \mathbb{E} \left [  1.15^T \right ]^{1-1/p}   1.15^{-k(1-1/p)} \\
  & = & K_{p,x}  \left ( 2  \left \{ \frac{2}{(\beta^\ast p + 1)(\beta^\ast p + 2)} \right\}^{1/p}  1.15^{1/p-1} \right )^k,
\end{eqnarray*}
where $K_{p,x}$ is a constant that only depends on $p$ and $x$ but on $k$.
Now, one can easily prove that for $p>1$ sufficiently close to 1, the term between brackets in the last display becomes strictly less than $1$.
Consequently, letting $\varepsilon > 0$ fixed, there exists an integer $k$ sufficiently large such that for every $t > 0$, 
\begin{eqnarray}  \label{premeqimpt}
t^{- \beta^\ast} 2^{k} \mathbb{E} \left [ N_{\overline{M}_{k} t -F_{k}}  (X_k ) \textbf{1}_{T > k} \right ] & \leq & \varepsilon. 
\end{eqnarray}

\subsection{Conclusion}

Observe that we have for every $t > 0$
 \begin{eqnarray*}
& &  t^{- \beta^\ast} 2^{k} \mathbb{E} \left [ \tilde{N}_{\overline{M}_{k} t -F_{k}}  (X_{k}) \textbf{1}_{T \leq k} \right ] \\
 & = & 2^{k} \mathbb{E} \left [   \textbf{1}_{T \leq k}  \mathbb{E} \left [ \left .    t^{- \beta^\ast}  \tilde{N}_{\overline{M}_{k} t -F_{k}}  (X_{k}) \right | \sigma(\overline{M}_{k}, F_{k}, T) \right ] \right ] \\
 &=& 2^{k} \mathbb{E} \left [   \textbf{1}_{T \leq k} (\overline{M}_{k}-t^{-1}F_{k})^{\beta^*}_{+} \mathbb{E} \left [ \left . (\overline{M}_{k}t-F_{k})^{-\beta^*}_{+} \tilde{N}_{\overline{M}_{k} t -F_{k}}  (X_{k}) \right | \sigma(\overline{M}_{k}, F_{k}, T) \right ] \right ],
\end{eqnarray*}
where $y_+$ denotes $y  \vee 0$. By Proposition~\ref{mc}, on the event $\{ T \leq k\}$, the r.v.\,\,$X_{k}$  is uniformly distributed over $[0,1]$ and independent of $M_1, \dots , M_{k}$ thus of $\overline{M}_{k}$. It is also independent of $F_{k}$ and $T$. Hence, letting $\theta$ be the map $s \mapsto 
\mathbb{E} [s_+^{-\beta^{\ast}} N_s(U)]$, where $U$ is a random variable uniformly distributed on $(0,1)$ independent of $N$, we  have:
 \begin{eqnarray*}
t^{- \beta^\ast} 2^{k} \mathbb{E} \left [ \tilde{N}_{\overline{M}_{k} t -F_{k}}  (X_{k}) \textbf{1}_{T \leq k} \right ] & = &  2^{k} \mathbb{E} \left [   \textbf{1}_{T \leq k} (\overline{M}_{k}-t^{-1}F_{k})^{\beta^*}_{+} \theta(  \overline{M}_{k} t -F_{k}  ) \right ].
\end{eqnarray*}
Applying Proposition~\ref{uni}, $\theta(\overline{M}_{k} t -F_{k})$ a.s.\,\,tends to a finite limit as $t \to \infty$. 
Hence  by dominated convergence
$ t^{- \beta^\ast} 2^{k} \mathbb{E} \left [ N_{\overline{M}_{k} t -F_{k}}  (X_{k}) \textbf{1}_{T \leq k} \right ]$  has a finite limit as $t \to \infty$. We deduce from this fact, (\ref{keyequation}) and  (\ref{premeqimpt}) that
 \begin{eqnarray*}
\limsup_{t \rightarrow \infty}  t^{- \beta^\ast} \mathbb{E} \left [ N_t(x)  \right ]  - \liminf_{t \rightarrow \infty}  t^{- \beta^\ast} \mathbb{E} \left [ N_t(x)  \right ]   & \leq & \varepsilon.
\end{eqnarray*}
Since that inequality holds for every $\varepsilon > 0$,  $t^{-\beta^\ast} \E[N_t(x)]$ has a finite limit as $t \to \infty$ which we denote by $n_{\infty}(x)$:
$$ n_{\infty}(x) = \lim_{t\to\infty} t^{-\beta^*} \E\big[N_{t}(x)\big].$$

\section{Identifying the limit} \label{sectioncinq}

In this section, we show that $ x  \mapsto n_\infty(x)$ is  proportional to  $x \mapsto (x(1-x))^{\beta^\ast/2}$ using a fixed point argument for integral equation (see also \cite[Section 4.1]{CLG10} for a similar application) . The normalizing constant will come from the $\mathbb{L}^1$-norm of $x\mapsto (x(1-x) )^{\beta^*/2}$ and the constant of Proposition~\ref{uni}. 

Combining (\ref{secoequaf}) with the densities computed in (\ref{loimx}), we deduce that
\begin{eqnarray*}
t^{-\beta^\ast}  \mathbb{E}[N_t(x)] & = &  t^{-\beta^\ast} \mathbb{P}(t \geq \tau_1)  + 2  \left ( \int_x^1 \frac{\textrm{d}y}{y} \int_0^{\frac{x}{y}} \textrm{d}m m^\bt \mathbb{E} \left [ (mt)^{-\bt}N_{m(t-\tau_1)}(y) \right ] \right. \\
& & + \left . \int_0^x \frac{\textrm{d}y}{1-y} \int_0^{\frac{1-x}{1-y}} \textrm{d}m m^\bt \mathbb{E} \left [(mt)^{-\bt}N_{m(t-\tau_1)}(y) \right ] \right ).
\end{eqnarray*}
Thanks to Lemma~\ref{lemmebornit}, we get by dominated convergence
 \begin{eqnarray*}
n_{\infty}(x) & = & \frac{2}{\bt+1} \left ( x^{\bt+1} \int_x^1 \textrm{d}y \frac{1}{y^{\bt+2}} n_{\infty}(y) + (1-x)^{\bt+1} \int_0^x \textrm{d}y \frac{1}{(1-y)^{\bt+2}} n_{\infty}(y)    \right ).
 \end{eqnarray*}
 In other words, if we define 
\begin{eqnarray*}
g_x(y) & = &  \frac{2}{\bt+1} \left ( x^{\bt+1} \frac{1}{y^{\bt+2}} \textbf{1}_{x < y < 1}  + (1-x)^{\bt+1} \frac{1}{(1-y)^{\bt+2}} \textbf{1}_{0 < y < x} \right ),
 \end{eqnarray*}
 we have
\begin{eqnarray*}
n_{\infty}(x) & = & \int_0^1 \textrm{d}y g_x(y) n_{\infty}(y).
 \end{eqnarray*}
Let $G$ be the operator that maps a function $f \in \mathbb{L}^1[0,1]$ to the function
 \begin{eqnarray*}
 G(f) (x) & = & \int_0^1 \textrm{d}y g_x(y) f(y).
 \end{eqnarray*}
 In particular, 
 $n_{\infty}$ is a fixed point of $G$. It is easy to check that $x\in (0,1)  \mapsto g_x(.) \in \mathbb{L}^1[0,1]$ is continuous for the $\mathbb{L}^1$-norm. Furthermore, Lemma~\ref{lemmebornit} ensures that $|n_\infty(x)| \leq C$ for every $x\in (0,1)$. As a consequence, $x\mapsto n_\infty(x)$ is continuous over $(0,1)$. 
An easy computation shows that for every $y \in (0,1)$, $\int_0^1\textrm{d}x g_x(y)=1$. 
Let $p$ be another fixed point of $G$ having the same integral as $n_\infty$. Then
\begin{eqnarray*}
\int_0^1 \textrm{d}x |n_\infty(x) - p(x)| & = &  \int_0^1 dx \left | \int_0^1 dy g_x(y)(n_\infty - p)(y) \right | \\
& \leq & \int_0^1 dx  \int_0^1 dy g_x(y) \left | n_\infty(y) - p(y) \right | \\ &=& \int_0^1 \textrm{d}y |n_\infty(y) - p(y)|,
\end{eqnarray*}
which shows that the inequality is in fact an equality. Hence $n_\infty - p$ has a.e.\,\,a constant sign. As we know that the integral  of $n_{\infty}-p$ is zero, we deduce that $n_\infty = p$ a.e. 
Straightforward calculations prove that $p_{0} : x \mapsto (x(1-x))^{\beta^\ast/2}$ is also a fixed point of $G$ of $\mathbb{L}^1$-norm, so that
\begin{eqnarray*}
n_{\infty}(x) &=& \|n_{\infty}\|_{1} \|p_{0}\|_{1}^{-1} \big(x(1-x)\big)^{\beta^\ast/2} \quad \textrm{a.e.}
\end{eqnarray*}
Since $n_\infty$ and $p_0$ are continuous, we can remove the a.e.\,\,statement (observe that $n_\infty(0)=n_\infty(1) = 0$ by Theorem~\ref{lastthro}).
Plainly,
$$
 \|p_{0}\|_{1}= \frac{\Gamma^2\left(\frac{\beta^*}{2}+1\right)}{\Gamma(\beta^*+2)}.
 $$
 On the other hand,  (\ref{bornitude}) and the dominated convergence theorem ensure that  $\|n_{\infty}\|_{1} = \lim_{t\to\infty}t^{-\beta^*}\E[ N_{t}(U)]$, which was computed in Proposition~\ref{uni}:
$$ \|n_{\infty}\|_{1}= \frac{\Gamma(2(\beta^*+1))}{2\Gamma^3(\beta^*+1)}.$$ 

\proof[Proof of Theorem~\ref{mainresul}] To sum up, we have for every $x \in [0,1]$: 
\begin{eqnarray*}
t^{-\beta^\ast} \mathbb{E}\big[N_{t}(x)\big] & \underset{t \to \infty}{\longrightarrow} & \frac{\Gamma\left(2\beta^*+2\right) \Gamma(\beta^*+2) }{ 2 \Gamma^3(\beta^*+1) \Gamma^2\left(\frac{\beta^*}{2}+1\right)} \big( x(1-x) \big) ^{\beta^*/2}.
 \end{eqnarray*}
Applying Lemma~\ref{relationnt2}, Theorem~\ref{mainresul} is shown.
\endproof


\section{Extensions and comments}

\subsection{Various convergences}

In this paper, we only proved a convergence in mean of $t^{-\beta^\ast} N_t(x)$. We may wonder whether this quantity also converges in distribution, in probability, or even almost surely. A more interesting question is the following: does the process $((t^{-\beta^\ast} N_t(x))_{x \in [0,1]}, t > 0)$ converge in distribution in the Skorokhod sense to a random function $(\mathcal{C}(x))_{x \in [0,1]}$ as $t \to \infty$? Observe that if it does, then there exists a random point $U$ uniformly distributed over $(0,1)$ such that $\mathcal{C}(U) = 0$, $U$ corresponding to the point $x_1$ of the first atom of $\Pi$ ($N_t(x_1)$ is indeed of order $t^{\sqrt{2}-1}$ by Theorem~\ref{lastthro}).

\begin{conjecture}
We have the functional limit law $(t^{-\beta^\ast} N_t(x))_{x \in [0,1]} \to (\mathcal{C}(x))_{x \in [0,1]}$ as $t \to \infty$  in $(\mathbb{D}([0,1]), \| \cdot \|_\infty)$, where $\mathcal{C}$ satisfies the distributional fixed point equation
\begin{eqnarray*}
 (\mathcal{C}(x))_{x \in [0,1]} & \overset{(d)}{=} & \left ( \mathbf{1}_{x < U_0} \left \{  \left ( U_0 U_1 \right )^{\beta^\ast} \mathcal{C}^{(00)} \left ( \frac{x}{U_0} \right ) +  \left ( U_0 (1- U_1) \right )^{\beta^\ast} \mathcal{C}^{(01)} \left ( \frac{x}{U_0} \right ) \right \}  \right.\\
& & + \mathbf{1}_{x > U_0} \left \{  \left ( (1-U_0) U_1 \right )^{\beta^\ast} \mathcal{C}^{(10)} \left ( \frac{x-U_0}{1-U_0} \right )  \right.  \\
& &  + \left. \left. \left ( (1-U_0) (1- U_1) \right )^{\beta^\ast} \mathcal{C}^{(11)} \left ( \frac{x-U_0}{1-U_0} \right ) \right \} \right )_{x \in [0,1]},
\end{eqnarray*}
where $U_0$, $U_1$, $\mathcal{C}^{(00)}$, $\mathcal{C}^{(01)}$, $\mathcal{C}^{(10)}$, $\mathcal{C}^{(11)}$ are independent, $U_0$ and $U_1$ are uniformly distributed on $[0,1]$ and $\mathcal{C}^{(00)}$, $\mathcal{C}^{(01)}$, $\mathcal{C}^{(10)}$, $\mathcal{C}^{(11)}$ have all the same distribution as $\mathcal{C}$.

\end{conjecture}

\subsection{Multidimensional case}

 The strategy adopted in Section~\ref{uniformsection} may be generalized to higher dimensions. As for the convergence in mean of the number of  hyper-rectangles crossed by a \emph{fixed} affine subspace having a direction generated by some vectors of the canonical basis, our approach may also be followed.
 \subsection{Quadtree as a model of random geometry}
 On top of its numerous applications in theoritical computer science, the model of random quadtree may be considered as a model of random geometry. More precisely one can view, for $t\geq 0$, the set of rectangles $\mathrm{Q}(t)$ as a random graph, assigning length $1$ to each edge of the rectangles. We denote this graph by $\tilde{\mathrm{Q}}(t)$. A natural question would be to understand the metric behavior of $\tilde{\mathrm{Q}}(t)$ as $t \to \infty$? If $L_{t}$ is the graph distance in $\tilde{\mathrm{Q}}(t)$ between the up-left and up-right corners, then Theorem \ref{lastthro} already shows that $L_{t}$ is less than the order $t^{\sqrt{2}-1}$. 
 \begin{problem} Is $\big(t^{1-\sqrt{2}}L_{t}\big)_{t\geq0}$ tight? Does it converge in distribution? If not, find the right power of $t$.
 \end{problem}

\addcontentsline{toc}{section}{References}
\bibliographystyle{abbrv}

\begin{thebibliography}{10}

\bibitem{Ber06}
J.~Bertoin.
\newblock {\em Random Fragmentations and Coagulation Processes}.
\newblock Number 102 in Cambridge Studies in Advanced Mathematics. Cambridge
  University Press, 2006.

\bibitem{BG04}
J.~Bertoin and A.~Gnedin.
\newblock Asymptotic laws for nonconservative self-similar fragmentations.
\newblock {\em Electron. J. Probab.}, 9(19):575--593, 2004.

\bibitem{CHH03}
H.-H. Chern and H.-K. Hwang.
\newblock Partial match queries in random quadtrees.
\newblock {\em SIAM J. Comput.}, 32(4):904--915 (electronic), 2003.

\bibitem{CLG10}
N.~Curien and J.-F. Le~Gall.
\newblock Random recursive triangulations of the disk \textit{via}
  fragmentation theory.
\newblock {\em preprint available on arxiv}, 2010.

\bibitem{FB74}
R.~A. Finkel and J.~L. Bentley.
\newblock Quad trees a data structure for retrieval on composite keys.
\newblock {\em Acta Informatica}, 4(1):1--9, mars 1974.

\bibitem{FGPR93}
P.~Flajolet, G.~Gonnet, C.~Puech, and J.~M. Robson.
\newblock Analytic variations on quadtrees.
\newblock {\em Algorithmica}, 10(6):473--500, 1993.

\bibitem{FLLS95}
P.~Flajolet, G.~Labelle, L.~Laforest, and B.~Salvy.
\newblock Hypergeometrics and the cost structure of quadtrees.
\newblock {\em Random Structures Algorithms}, 7(2):117--144, 1995.

\bibitem{FS09}
P.~Flajolet and R.~Sedgewick.
\newblock {\em Analytic combinatorics}.
\newblock Cambridge University Press, Cambridge, 2009.

\bibitem{IL97}
A.~Iserles and Y.~Liu.
\newblock Integro-differential equations and generalized hypergeometric
  functions.
\newblock {\em J. Math. Anal. Appl.}, 208(2):404--424, 1997.

\bibitem{Liu97}
Q.~Liu.
\newblock Sur une \'equation fonctionnelle et ses applications: une extension
  du th\'eor\`eme de {K}esten-{S}tigum concernant des processus de branchement.
\newblock {\em Adv. in Appl. Probab.}, 29(2):353--373, 1997.

\bibitem{Liu01}
Q.~Liu.
\newblock Asymptotic properties and absolute continuity of laws stable by
  random weighted mean.
\newblock {\em Stochastic Process. Appl.}, 95(1):83--107, 2001.

\bibitem{MT09}
S.~Meyn and R.~L. Tweedie.
\newblock {\em Markov chains and stochastic stability}.
\newblock Cambridge University Press, Cambridge, second edition, 2009.
\newblock With a prologue by Peter W. Glynn.

\end{thebibliography}

\end{document}